# Accelerating Fleet Upgrade Decisions with Machine-Learning Enhanced Optimization

Kenrick Howin Chai[a,*], Stefan Hildebrand[a,*], Tobias Lachnit[b], Martin Benfer[b], Gisela Lanza[b], Sandra Klinge[a]

[a]Department of Structural and Computational Mechanics, TU Berlin, Straße des 17. Juni 135, Berlin and 10623, Germany
[b]wbk Institute of Production Science, Karlsruhe Institute of Technology, Kaiserstraße 12, Karlsruhe and 76131, Germany

* Corresponding author. Tel.: +49 (0) 30 314 21481; E-mail address: stefan.hildebrand@tu-berlin.de, kenrick.howin.chai@campus.tu-berlin.de

**Abstract**

Rental-based business models and increasing sustainability requirements intensify the need for efficient strategies to manage large machine and vehicle fleet renewal and upgrades. Optimized fleet upgrade strategies maximize overall utility, cost, and sustainability. However, conventional fleet optimization does not account for upgrade options and is based on integer programming with exponential runtime scaling, which leads to substantial computational cost when dealing with large fleets and repeated decision-making processes. This contribution firstly suggests an extended integer programming approach that determines optimal renewal and upgrade decisions. The computational burden is addressed by a second, alternative machine learning-based method that transforms the task to a mixed discrete-continuous optimization problem. Both approaches are evaluated in a real-world automotive industry case study, which shows that the machine learning approach achieves near-optimal solutions with significant improvements in the scalability and overall computational performance, thus making it a practical alternative for large-scale fleet management.

## 1. Introduction

With rising sustainability demands and the need for long-term cost-effective operations, extending the lifetime of capital equipment has become increasingly important. Rather than relying on costly and recurring replacements, companies are increasingly exploring individualized upgrades and adaptations to maintain or enhance fleet performance over time. In Product-Service Systems (PSS), enabling mid-life upgrades provides a promising strategy to reduce resource consumption, improve operational efficiency, and ensure continuous technological adaptation [1]. This shift in lifecycle management is particularly relevant in the context of Industry 4.0, where rapidly evolving digital technologies must be integrated into long-lived mechanical systems [2]. In fleet management, decisions still often focus on full replacements, although targeted upgrades represent a valuable alternative or addition to conventional renewal strategies. However, planning such upgrades at the fleet level introduces substantial complexity for manufacturers. Varying usage patterns, system aging, and heterogeneous product states must be considered while maintaining economic viability and scalability. Due to the scale of product fleets and the multitude of upgrade possibilities, conventional optimisation approaches like mixed-integer programming (MIP) often face significant computational challenges. To address this challenge, the present study explores how machine learning (ML) can accelerate decision-making in fleet renewal scenarios by enabling the efficient planning and control of individualized product upgrades and adaptations at scale.

The conventional Fleet Renewal Problem (FRP), also referred to as fleet replacement problem [3] asks for the optimal timing and selection of vehicle or machine replacements and has been studied extensively [4–7].

Extended approaches focus to capture uncertainty in key parameters [8–10], to include varying usage intensity, technological changes and asset deterioration over the product lifespan [11,12] and environmental considerations embedded into the objective function or modeled as constraints [13–15].

Overall, the existing research on combined fleet renewal and upgrade decisions is sparse as predominantly upgrades are addressed at the level of individual products, with limited attention to fleet-level considerations [16–19]. This results in fragmented strategies and suboptimal outcomes which is particularly limiting in scenarios where both replacement and upgrade options are viable alternatives.

A major challenge when including product upgrades in management strategies for large fleets is the computational effort as FRP is solved conventionally discretizing all variables including the age of the products into integer values, leading to a combinatoric problem in the complexity class NP. A wide range of optimization methods has been proposed, most commonly integer programming (IP) [11,15] and MIP [12], where the computational effort grows exponentially with the problem size, intensified by intense RAM usage.

The present study addresses the direct inclusion of upgrade strategies into FRP, respecting both cost and sustainability objectives, by a unified ML approach based on Straight-Through Estimator (STE). The new approach allows for significantly improved scalability using gradient-based optimization with only polynomial computational complexity. The validation is carried out with a real-world automotive industry case study. Its results show that the ML approach delivers near-optimal solutions with significantly improved computational performance, making it a practical tool for modern fleet lifecycle management.

## 2. Conceptual und Methodological Background

A product upgrade is a strategic approach to extend a product's lifecycle by enhancing its functionality and performance to meet evolving requirements, such as technological advancements or changing customer preferences. Unlike complete product replacement, upgrades focus on modifying or improving specific modules or components, allowing for continued use of the existing system [16,20]. Existing research covers financial aspects [1], design aspects [21–23], maintenance strategies for dealers and the second-hand market [24,25] as well as upgrade strategies for end-users [16,26].

The term fleet broadly refers to a collection of product entities that are grouped based on shared characteristics or a common operational purpose. While traditionally associated with ships or vehicles, the concept extends to encompass systems, sub-systems, or components within an industrial framework, depending on the focus of analysis. The defining feature of a fleet is the commonality among its members, whether technical, operational or contextual [27].

The FRP asks for an optimal schedule to replace units within a fleet which minimizes the total cost of ownership and operation over a defined planning horizon. As vehicles age or accumulate usage, their performance deteriorates, leading to increased expenses such as operating, maintenance, repair and fixed overhead costs. The objective of the FRP is to strategically plan replacements in a way that balances these rising costs with investment in new assets, ensuring overall cost-efficiency and operational effectiveness of the fleet [11].

A point $(x^* \in F)$ in the feasible region $F \subset G \subset R^n, n \in N$ is considered optimal if it fulfills

$$x \in F \Rightarrow \varphi(x^*) \leq \varphi(x) \tag{1}$$

for an objective function $\varphi: G \to R$ [28].

IP and MIP problems aim to optimize a linear objective function

$$\varphi(x) := \sum_{i=1}^{n} \gamma_i \xi_i, \quad x = (\xi_1, \dots, \xi_n)^T \in R^n, \tag{2}$$

subject to a set of linear constraints

$$F := \{x \in R^n : a_i^T x \leq \beta_i, i = 1, \dots, m\}, \tag{3}$$

where $\gamma_i \in R$ are the objective coefficients, $a_i \in R^n$ the constraint vectors, and $\beta_i \in R$ the constraint bounds.

The distinction between IP and MIP lies in the domain of the decision variables

$$\xi_j \in \begin{cases} Z \text{ for all } j \in \{1, \dots, n\} & \text{(IP)} \\ Z \text{ for all } j \in J \subseteq \{1, \dots, n\} & \text{(MIP)} \end{cases}, \tag{4}$$

i.e. in IP all decision variables are restricted to integer values, while in MIP only a subset $J$ of variables is integer-constrained and the remaining variables can assume continuous values. Additional constraints such as equalities $a_i^T x = \beta_i$, inequalities $a_j^T x \geq \beta_j$ and variable bounds like $\xi_k \geq 0$ may also be included.

In the following, the framework to solve the conventional FRP is introduced [11,29]. IP Optimization is conducted over a planning horizon $T$, with replacement decisions made at each period $j \in \{0, \dots, T\}$. Each asset is characterized by its type $o \in \{0, \dots, O\}$ and age $i \in \{0, \dots, N\}$, where $N$ denotes the maximum useful life and the following further notations are:

*Decision Variables*

- $v_{oij}$    Number of type $o$, age $i$ assets deployed in period $j$.
- $b_{oij}$    Number of type $o$, age $i$ assets purchased in period $j$.
- $s_{oij}$    Number of type $o$, age $i$ assets sold in period $j$.
- $\delta_j$    Binary variable, which takes a value of 1 if any asset is purchased in period $j$ or 0 otherwise.

*Parameters*

- $f$    Discount factor per period.
- $p_{oij}$    Purchase price of a type $o$, age $i$ assets in period $j$.
- $k_j$    Fixed cost incurred in period $j$ if a purchase is made.
- $c_{oij}$    Operation and maintenance cost of a type $o$, age $i$ asset in period $j$.

$r_{oij}$    Resale or salvage value of a type $o$, age $i$ asset in period $j$.
$u_i$    Usage capacity of an asset of age $i$ per period.
$d_j$    Fleet capacity requirement in period $j$.
$h_{oi}$    Initial number of type $o$, age $i$ assets at the start of the planning horizon.

*Objective Function*

$$\text{minimize} \quad C_{\text{objective,base}} = \sum_{j=0}^{T-1} f^j \left( \sum_{i=0}^{N-1} \sum_{o=0}^{O} p_{oij} b_{oij} + k_j \delta_j \right)$$
$$+ \sum_{j=0}^{T-1} \sum_{i=0}^{N-1} \sum_{o=0}^{O} f^{j+1} c_{oij} v_{oij} - \sum_{j=0}^{T} \sum_{i=0}^{N} \sum_{o=0}^{O} f^j r_{oij} s_{oij}$$

subject to

$$\sum_{i=0}^{N-1} \sum_{o=0}^{O} u_i v_{oij} \geq d_j \qquad \forall j < T \qquad (5)$$
$$b_{oij} - v_{oij} = 0 \qquad \forall o; i = 0; \forall j < T \qquad (6)$$
$$h_{oi} + b_{oij} = v_{oij} + s_{oij} \qquad \forall o; \forall 0 < i < N; j = 0 \qquad (7)$$
$$v_{o,i-1,j-1} + b_{oij} = v_{oij} + s_{oij} \qquad \forall o; \forall 0 < i < N; \forall 0 < j < T \qquad (8)$$
$$h_{oi} - s_{oij} = 0 \qquad \forall o; i = N; j = 0 \qquad (9)$$
$$v_{o,i-1,j-1} - s_{oij} = 0 \qquad \{\forall o; i = N; \forall j\} \text{ or } \{\forall o; \forall i > N; j = T\} \qquad (10)$$
$$b_{o,i,j} = 0 \qquad \forall o; \forall i; j = T \qquad (11)$$
$$s_{o,i,j} = 0 \qquad \forall o; i = 0; \forall j \qquad (12)$$
$$\sum_{i=0}^{N-1} \sum_{o=0}^{O} b_{oij} \leq M \delta_j \qquad \forall j < T \qquad (13)$$
$$v_{oij}, b_{oij}, s_{oij} \in \mathbb{N}_{\geq 0} \qquad \forall o; \forall i; \forall j \qquad (14)$$
$$\delta_j \in \{0,1\} \qquad \forall o; \forall i; \forall j \qquad (15)$$

Constraints (5) ensure that the fleet meets the required capacity in every period. Constraint (6) enforces that newly acquired assets are immediately deployed. The flow balance for the initial period is modeled in constraint (7), while constraint (8) governs the asset transitions across subsequent periods. Constraint (9) ensures that initial assets reaching their maximum life are sold in the first period. Constraint (10) guarantees that all assets are sold by the end of their lifetime or at the end of the planning horizon. Constraint (11) prohibits purchases in the final period, and constraint (12) prevents immediate resale of newly acquired assets. Constraint (13) apply fixed purchasing costs to each period in which at least one asset is purchased. Finally, constraints (14) and (15) define the domains of the decision variables.

## 3. Fleet Upgrade Problem for sustainability

Fleet upgrading is particularly beneficial for high-investment, long-lifetime products, as it reduces costs and environmental impact associated with full disposal and repurchase. While maintenance can slow functional decline, targeted upgrades can not only preserve but enhance product capabilities beyond the original design, offering a cost-effective and sustainable alternative to traditional replacement cycles [16].

In the following, the framework of the FRP is extended to model product upgrades besides product replacements in a fleet and to consider environmental costs in the objective function. While the general structure of the model remains unchanged, several modifications are introduced. Specifically, each possible configuration of a product is treated as a distinct type, denoted by $o \in O$. For example, $o = 0$ represents the base configuration, $o = 1$ the first upgraded variant, $o = 2$ the second upgrade, and so on. To accommodate this extension, additional decision variables and parameters are introduced, and the existing constraints are either modified or supplemented with new ones to accurately reflect the upgrade dynamics.

*Decision Variables*

$m_{oij}$    Number of type $o$ components installed on age $i$ assets in period $j$.
$dm_{oij}$    Number of type $o$ components removed from age $i$ assets in period $j$.
$bk_{oij}$    Number of new components of type $o$ purchased in period $j$.
$sk_j$    Number of used components of type $o$ sold in period $j$.
$i_{oj}$    Inventory of type $o$ components held in period $j$.

*Parameters*

$ep_{oj}$    Emissions from producing a type $o$ asset in period $j$ (tons $CO_2$-eq).
$en_{oj}$    Emissions during operation of a type $o$ asset in period $j$ (tons $CO_2$-eq).
$es_{oj}$    Emissions from disposal of a type $o$ asset in period $j$ (tons $CO_2$-eq).
$ec_j$    $CO_2$ price in period $j$ (cost per ton $CO_2$-eq).
$hk_o$    Initial inventory of type $o$ components before planning starts.
$pk_{oj}$    Purchase price of a type $o$ component in period $j$.
$rk_{oj}$    Selling price of a type $o$ component in period $j$.
$cm_{oj}$    Installation cost of a type $o$ component in period $j$.
$cd_{oj}$    Removal cost of a type $o$ component in period $j$.
$u_{oij}$    Usage capacity of a type $o$, age $i$ asset in period $j$.

*Objective Function*

$$\text{minimize} \quad C_{\text{objective,extended}} = C_{\text{objective,base}} + C_{\text{environment}} + C_{\text{upgrade}}$$

$$C_{\text{environment}} = \sum_{j=0}^{T-1} \sum_{i=0}^{N-1} \sum_{o=0}^{O} f^j \cdot ep_{oj} \cdot ec_j \cdot b_{oij}$$
$$+ \sum_{j=0}^{T-1} \sum_{i=0}^{N-1} \sum_{o=0}^{O} f^{j+1} \cdot en_{oj} \cdot ec_j \cdot v_{oij} + \sum_{j=0}^{T} \sum_{i=0}^{N} \sum_{o=0}^{O} f^j \cdot es_{oj} \cdot ec_j \cdot s_{oij}$$

$$C_{\text{upgrade}} = \sum_{o=0}^{O} \sum_{j=0}^{T} f^j \cdot pk_{oj} \cdot bk_{oj} + \sum_{o=0}^{O} \sum_{i=0}^{N} \sum_{j=0}^{T} f^j \cdot cm_{oj} \cdot m_{oij}$$
$$+ \sum_{o=0}^{O} \sum_{i=0}^{N} \sum_{j=0}^{T} f^j \cdot cd_{oj} \cdot dm_{oij} - \sum_{o=0}^{O} \sum_{j=0}^{T} f^j \cdot rk_{oj} \cdot sk_{oj}$$

subject to

$$\sum_{i=0}^{N-1} \sum_{o=0}^{O} u_i v_{oij} \geq d_j \qquad \forall j < T \qquad (16)$$
$$\sum_{o=0}^{O} (b_{oij} - v_{oij}) = 0 \qquad i = 0; \forall j < T \qquad (17)$$

$$\sum_{o=0}^{O}(h_{oi} + b_{oij}) = \sum_{o=0}^{O}(v_{oij} + s_{oij}) \qquad \forall i; j = 0 \qquad (18)$$

$$\sum_{o=0}^{O} b_{oij} = \sum_{o=0}^{O}(v_{oij} + s_{oij}) \qquad i = 0; \forall 0 < j \le T \qquad (19)$$

$$\sum_{o=0}^{O}(v_{o,i-1,j-1} + b_{oij}) = \sum_{o=0}^{O}(v_{oij} + s_{oij}) \qquad \forall 0 < i \le N; \forall 0 < j \le T \qquad (20)$$

$$m_{oij} + b_{oij} = v_{oij} - h_{oi} + dm_{oij} + s_{oij} \qquad \forall o; \forall i; j = 0 \qquad (21)$$

$$m_{oij} + b_{oij} = v_{oij} + dm_{oij} + s_{oij} \qquad \forall o; i = 0; \forall 0 < j \le T \qquad (22)$$

$$m_{oij} + b_{oij} = v_{oij} - v_{o,i-1,j-1} + dm_{oij} + s_{oij} \qquad \forall o; \forall 0 < i \le N; \forall 0 < j \le T \qquad (23)$$

$$i_{oj} = hk_o + bk_{oj} - sk_{oj} + \sum_{i=0}^{N}(dm_{oij} - m_{oij}) \qquad \forall o; j = 0 \qquad (24)$$

$$i_{oj} = i_{o,j-1} + bk_{oj} - sk_{oj} + \sum_{i=0}^{N}(dm_{oij} - m_{oij}) \qquad \forall o; 0 < j \le T \qquad (25)$$

$$\sum_{o=0}^{O}(h_{oi} - s_{oij}) = 0 \qquad i = N; j = 0 \qquad (26)$$

$$\sum_{o=0}^{O}(v_{o,i-1,j-1} - s_{oij}) = 0 \qquad \{i = N; \forall j\} \text{ or } \{\forall i > N; j = T\} \qquad (27)$$

$$b_{o,i,j} = 0 \qquad \forall o; \forall i; j = T \qquad (28)$$
$$s_{o,i,j} = 0 \qquad \forall o; i = 0; \forall j \qquad (29)$$
$$b_{o,i,j} = 0 \qquad \forall o \ne 0; \forall 0 < i \le N; \forall j \qquad (30)$$
$$s_{o,i,j} = 0 \qquad \forall o \ne 0; \forall 0 < i \le N; \forall j \qquad (31)$$

$$\sum_{i=0}^{N-1}\sum_{o=0}^{O} b_{oij} \le M\delta_j \qquad \forall j < T \qquad (32)$$

$$v_{oij}, b_{oij}, s_{oij}, m_{oij}, dm_{oij}, bk_{oj}, sk_{oj}, I_{oj} \in \mathbb{N}_{\ge 0} \quad \forall o; \forall i; \forall j \qquad (33)$$
$$\delta_j \in \{0,1\} \qquad \forall o; \forall i; \forall j \qquad (34)$$

Constraints (16) ensure that the fleet satisfies the required capacity in each planning period. Constraint (17) enforces that newly acquired assets are immediately deployed. The initial asset flow balance is captured by constraint (18), while the dynamic flow of assets across subsequent periods is governed by constraints (19) and (20). Constraints (21) to (23) model the upgrade dynamics, covering both the initial configuration and period-to-period transitions. Component stock balance is established in constraint (24), with ongoing inventory adjustments managed through constraint (25). Constraint (26) mandates the disposal of initial assets that have reached their maximum service life. Constraint (27) ensures that all assets are sold either upon reaching their maximum age or by the end of the planning horizon. Constraint (28) prohibits any asset purchases in the final period, while constraint (29) prevents immediate resale of newly acquired assets, ensuring at least one period of use. Optional constraints (30) and (31) restrict asset transactions to the base configuration, should such a policy be desired. Constraint (32) apply fixed purchasing costs to each period in which at least one asset is purchased. Finally, constraints (33) and (34) specify the admissible domains for all decision variables, including integrality and binary requirements.

## 4. Accelerating the Optimization with Machine Learning techniques

IP is well known for being NP-hard [30], i.e. the required computational complexity grows exponentially with the problem size. To mitigate this computational burden, an alternative approach is suggested which reformulates the IP into a mixed discrete-continuous optimization problem. This reformulated problem is then solved with established ML techniques. Therefore, the STE is applied which introduces a differentiable approximation for discrete variables to enable the use of gradient based optimization algorithms with improved convergence behavior [31,32].

The extended model taking product upgrades into account entails disjunctive or logical constraints which cannot be directly included with standard IP solvers but requires explicit modeling techniques such as Big-M or hull transformation [33]. In preliminary tests, these increased the runtimes of the solver about more than three orders of magnitude even for small tasks, rendering this approach practically infeasible. In contrast, the ML framework presented here allows to include these directly. Leveraging this flexibility, the present formulation treats the asset allocation $v$ as the central decision variable to be directly optimized. All other decision variables are not optimized independently but are instead derived from $v$. This way, the original flow balance constraints can be omitted, which govern variable transitions across time periods in the IP formulation. This transformation reduces the complexity of the feasible set and eliminates potential multicollinearity among the optimized decision variables, thereby enhancing the efficiency and numerical stability of the optimization process.

Constraints that do not govern flow balance are reformulated as soft constraints and penalties for violations are calculated using the Rectified Linear Unit (ReLU) function, resulting in continuous functions facilitating the application of gradient-based optimizers.

The following equations define the mapping from the central variable $v$ to the full set of decision variables for the basic integer programming formulation

$$b_{oij} = \max(0, v_{oij} - h_{oi}) \qquad \forall o; \forall i; j = 0 \qquad (35)$$
$$b_{oij} = \max(0, v_{oij}) \qquad \forall o; \forall i = 0; \forall 0 < j \le T \qquad (36)$$
$$b_{oij} = \max(0, v_{oij} - v_{o,i-1,j-1}) \qquad \forall o; \forall 0 < i \le N; \forall 0 < j \le T \qquad (37)$$
$$s_{oij} = \max(0, h_{oi} - v_{oij}) \qquad \forall o; \forall i; j = 0 \qquad (38)$$
$$s_{oij} = \max(0, v_{o,i-1,j-1} - v_{oij}) \qquad \forall o; \forall 0 < i \le N; \forall 0 < j \le T \qquad (39)$$

$$\delta_j = \min\left(1, \sum_{o=0}^{O}\sum_{i=0}^{N-1} b_{oij}\right) \qquad \forall j < T \qquad (40)$$

For the extended integer programming formulation, the decision variables are derived from $v$ by the following equations

$$b_{oij} = \max\left(0, \sum_{o=0}^{O}(v_{oij} - h_{oi})\right) \qquad o = 0; \forall i; j = 0 \qquad (41)$$

$$b_{oij} = \max\left(0, \sum_{o=0}^{O} v_{oij}\right) \qquad \begin{aligned}&o = 0; i = 0;\\ &\forall 0 < j \le T\end{aligned} \qquad (42)$$

$$b_{oij} = \max\left(0, \sum_{o=0}^{O}(v_{oij} - v_{o,i-1,j-1})\right) \qquad \begin{aligned}&o = 0; \forall 0 < i \le N;\\ &\forall 0 < j \le T\end{aligned} \qquad (43)$$

$$s_{oij} = \max\left(0, \sum_{o=0}^{O}(h_{oi} - v_{oij})\right) \qquad o = 0; \forall i; j = 0 \qquad (44)$$

$$s_{oij} = \max\left(0, \sum_{o=0}^{O}(v_{o,i-1,j-1} - v_{oij})\right) \qquad \begin{aligned}&o = 0; \forall 0 < i \le N;\\ &\forall 0 < j \le T\end{aligned} \qquad (45)$$

$$m_{oij} = \max(0, v_{oij} - h_{oi} + s_{oij} - b_{oij}) \qquad \forall o; \forall i; j = 0 \qquad (46)$$
$$m_{oij} = \max(0, v_{oij} + s_{oij} - b_{oij}) \qquad \forall o; i = 0; \forall 0 < j \le T \qquad (47)$$
$$m_{oij} = \max(0, v_{oij} - v_{o,i-1,j-1} + s_{oij} - b_{oij}) \qquad \begin{aligned}&\forall o; \forall 0 < i \le N;\\ &\forall 0 < j \le T\end{aligned} \qquad (48)$$

$$d_{oij} = \max(0, h_{oi} - v_{oij} + b_{oij} - s_{oij}) \qquad \forall o; \forall i; j = 0 \qquad (49)$$
$$d_{oij} = \max(0, -v_{oij} + b_{oij} - s_{oij}) \qquad \forall o; i = 0; \forall 0 < j \le T \qquad (50)$$

$$d_{oij} = \max(0, v_{o,i-1,j-1} - v_{oij} + b_{oij} - s_{oij}) \qquad \begin{aligned}&\forall o; \forall 0 < i \le N;\\ &\forall 0 < j \le T\end{aligned} \qquad (51)$$

$$bk_{oj} = \max\left(0, \sum_{i=0}^{N}(m_{oij} - d_{oij}) - hk_o\right) \qquad \forall o; j = 0 \qquad (52)$$

$$bk_{oj} = \max\left(0, \sum_{i=0}^{N}(m_{oij} - d_{oij}) - i_{o,j-1}\right) \quad \forall o; \forall 0 < j \leq T \quad (53)$$

$$i_{oj} = \sum_{i=0}^{N}(d_{oij} - m_{oij}) + hk_o + bk_{oj} \quad \forall o; j = 0 \quad (54)$$

$$i_{oj} = \sum_{i=0}^{N}(d_{oij} - m_{oij}) + i_{o,j-1} + bk_{oj} \quad \forall o; \forall 0 < j < T \quad (55)$$

$$i_{oj} = 0 \quad \forall o; j = T \quad (56)$$

$$sk_{oj} = \sum_{i=0}^{N}(d_{oij} - m_{oij}) + i_{o,j-1} \quad \forall o; j = T \quad (57)$$

$$\delta_j = \min\left(1, \sum_{o=0}^{O}\sum_{i=0}^{N-1} b_{oij}\right) \quad \forall j < T \quad (58)$$

The optimization process is structured into four sequential phases, each aimed at iteratively improving the central decision variable $v$ to find an optimal solution within the feasible region:

1. **Initialization**
   The process begins with the initialization of the central decision variable $v$, usually randomly sampled and as far as possible based on problem-specific heuristics. A suitable optimization algorithm such as Adam is selected to iteratively update $v$. Adam is a widely used adaptive stochastic optimization method. It extends the stochastic gradient descent optimizer by computing individual learning rates from estimating first and second moments of the gradients [34].

2. **Projection and Derivation of Decision Variables**
   To facilitate optimization, $v$ is first transformed into a feasible space by eliminating negative values and ensuring no assets are active in the final time period or at their maximum lifetime:

$$\begin{aligned}&v = \max(0, v)\\&v_{oij} = 0 \quad \forall o; \forall i; j = T,\\&v_{oij} = 0 \quad \forall o; i = N; \forall j\end{aligned} \quad (59)$$

   Next, to allow gradient flow while rounding values down to the nearest integer, $v$ is discretized using STE incorporating the definition from Fan et al. [35] $STE_{round}(v) = $ stop_grad($\lfloor v \rfloor$) − stop_grad($v$) + $v$.
   The operator `stop_grad` blocks gradient flow during backpropagation, so that the derivative is considered as 0. In the forward pass, this yields $STE_{round}(v) = \lfloor v \rfloor$. During backpropagation, gradients are computed using the surrogate $\nabla STE_{round}(v) = \nabla v$.
   After projection and discretization, all remaining decision variables are derived from both the continuous and rounded versions of $v$ using predefined equations.

3. **Cost Function and Penalty Term Calculation**
   The cost function is computed for both the fractional and rounded solutions. Additionally, a penalty term is added for the rounded solution to account for any constraint violations not handled during projection or derivation,
   for inequality constraints: $p_{ineq} = c \cdot \text{ReLU}(z)$,
   for equality constraints: $p_{eq} = c \cdot |z|$ with the penalty weight $c$, and the respective expression $z$.

   If the penalty is zero and the rounded solution yields the lowest objective value so far, it is logged as the current best solution.

4. **Backpropagation**
   The loss function combining the cost of the fractional solution and the penalty from the rounded solution is computed together with its gradients w.r.t. the optimizable parameters. One optimization step is carried out by the optimizer to update $v$. The process then loops back to Phase 2 for the next iteration.

## 5. Case Studies

### 5.1. Case 1: Basic Formulation

To allow for direct comparison and benchmarking, the first case study applies the data from Parthanadee et al. [11] considering a vehicle fleet across eight scenarios, varying in vehicle prices, operation and maintenance (O&M) costs, and vehicle utilization. Vehicle prices follow either constant depreciation or real price trends from used-car markets. O&M costs are either fixed or based on data from 177 passenger cars. Two usage models are considered: one with constant mileage, the other based on user preferences. Table 1 outlines the scenario structure.

Table 1. Overview of the scenarios.

| Sce. | Vehicle Price | O&M | Utilization | Demand (km) |
|---|---|---|---|---|
| 1 | Linear | Constant | Constant | 300.000 |
| 2 | Linear | Constant | User Preference | 228.000 |
| 3 | Linear | Empirical | Constant | 300.000 |
| 4 | Linear | Empirical | User Preference | 228.000 |
| 5 | Empirical | Constant | Constant | 300.000 |
| 6 | Empirical | Constant | User Preference | 228.000 |
| 7 | Empirical | Empirical | Constant | 300.000 |
| 8 | Empirical | Empirical | User Preference | 228.000 |

The planning horizon is three years, with a maximum vehicle age of ten years and a 5% discount rate. The initial fleet consists of vehicles aged 1, 3, 5, 7, and 9 years, with 3 vehicles per age group. Each vehicle can provide 20,000 km of services In scenarios with constant mileage, total annual demand is 300,000 km; in preference-based scenarios, 228,000 km. Table 2 summarizes vehicle prices (in THB), O&M costs (THB/km), and lifetime mileage.

Table 2. Vehicle Price, O&M-Cost und Utilization.

| Age | Vehicle Price (THB) | | O&M-Cost (THB/km) | | Utilization (km) | |
|---|---|---|---|---|---|---|
| | Linear | Empirical | Constant | Empirical | Constant | User Preference |
| 0 | 1.000.000 | 1.000.000 | 3 | 3,44 | 20.000 | 20.000 |
| 1 | 946.315 | 910.000 | 3 | 3,53 | 20.000 | 20.000 |
| 2 | 892.630 | 830.000 | 3 | 4,75 | 20.000 | 20.000 |
| 3 | 838.945 | 820.000 | 3 | 3,75 | 20.000 | 20.000 |
| 4 | 785.260 | 630.000 | 3 | 4,81 | 20.000 | 18.000 |
| 5 | 731.575 | 620.000 | 3 | 3,82 | 20.000 | 16.000 |
| 6 | 677.890 | 600.000 | 3 | 3,73 | 20.000 | 14.000 |
| 7 | 624.205 | 590.000 | 3 | 4,91 | 20.000 | 12.000 |
| 8 | 570.520 | 580.000 | 3 | 3,83 | 20.000 | 10.000 |
| 9 | 516.835 | 510.000 | 3 | 5,02 | 20.000 | 8.000 |
| 10 | 463.150 | 330.000 | 3 | 4,05 | 0 | 0 |

*5.2. Case 2: Extended Formulation*

An electric vehicle with exchangeable battery modules is introduced into the second case study to represent an upgrade option. The chosen example product is the Nio ES8 which is currently available in the market with two battery configurations: a 75 kWh battery providing a range of 375 km, and a 100 kWh battery with a range of 500 km. The base vehicle costs €60,000, excluding the battery. The battery costs amount to €12,000 for the 75 kWh version and €21,000 for the 100 kWh version [36]. In addition, each battery swap incurs a cost of €30 [37]. Annual maintenance costs are assumed to be €576. The annual mileage per vehicle is set at 11,733 km.

Depreciation is modeled as follows: the vehicle loses 25% of its value in the first years. After three years, 50% of the original value remains. From the fourth year onward, a constant annual depreciation rate of 5% is assumed. The parameters for environmental costs are summarized in Table 3, based on Wietschel et al. [38]. Emissions from disposal are neglected, as their contribution is negligible compared to other phases of the vehicle life cycle. The assumed $CO_2$ prices are €45 per tonne in 2024, €55 per tonne in 2025, and €65 per tonne in 2026 [39].

Table 3. $CO_{2eq}$-Emission from diesel and electric cars

| Type | 2024 | 2025 | 2026 |
|---|---|---|---|
| Diesel car production (t $CO_{2eq}$ per car) | 8,1 | 8,1 | 8,1 |
| Diesel car operation (g $CO_{2eq}$/km) | 134,68 | 136,36 | 138,04 |
| Electric car production (t $CO_{2eq}$ per car) | 12,6 | 12,5 | 12,4 |
| Electric car operation (g $CO_{2eq}$/km) | 5,12 | 5,05 | 4,97 |
| $CO_2$ price (€ / t) | 45 | 55 | 65 |

The planning horizon spans three years, with a maximum vehicle service life of ten years. At the start of the planning period, the fleet consists of no vehicles. The initial range requirement of the fleet is 5,000 km, increasing by 200 km annually.

# 6. Results

The optimization is carried out with CPLEX as representative of the IP method and with Adam for the ML method. The discrepancy is evaluated as absolute percentage error

$$\text{Discrepancy} = \frac{|\text{obj} - \text{obj}_{\text{optimal}}|}{|\text{obj}_{\text{optimal}}|} \cdot 100\% \quad . \quad (60)$$

The optimization results are presented in Table 4. The objective values for the first base case are negative, as the initial fleet consists of 15 vehicles and the investment costs are excluded from the objective function. The results indicate that the ML approach achieved the same optimal value as IP in 5 out of 9 scenarios. In the remaining cases, the ML method produced near-optimal solutions, with a maximum discrepancy of approximately 1%.

Furthermore, an empirical analysis is performed to investigate how computational complexity scales with problem size. For this purpose, scenario 7 from the base model is employed as an exemplary case study. The problem size is scaled by extending the planning horizon following a logarithmic progression, allowing a systematic investigation of the scalability of the optimization methods across multiple orders of magnitude. The discount rate is omitted in this analysis, as its application would lead to cost values approaching zero over extended planning periods, thereby distorting the results.

Table 4. Optimization result.

| Case | Sce. | Integer Programming | Machine Learning | Discrepancy |
|---|---|---|---|---|
| Base | 1 | -5.274.136 | -5.274.136 | 0,00% |
|  | 2 | -6.029.618 | -5.990.568 | 0,65% |
|  | 3 | -4.486.399 | -4.486.399 | 0,00% |
|  | 4 | -5.565.468 | -5.506.788 | 1,05% |
|  | 5 | -6.285.552 | -6.285.552 | 0,00% |
|  | 6 | -7.023.553 | -7.006.278 | 0,002% |
|  | 7 | -5.668.737 | -5.668.737 | 0,00% |
|  | 8 | -6.234.137 | -6.223.408 | 0,17% |
| Extended | - | 257.454 | 257.454 | 0,00% |

Fig. 1 compares the computation times of the IP solver and the ML solver depending on problem size. The results demonstrate that the IP solver handles small problems very rapidly. However, its computation time increases rapidly with the problem size as expected due to IP being an NP hard problem.

In contrast, although the ML method requires more time for smaller problems, its computational time grows at a significantly lower rate. This corresponds to the upper bound of training complexity with Adam not being significantly dependent on the problem size but being dominated by the required final tolerance [40]. Furthermore, the IP solver was unable to provide a solution starting at a planning horizon of $10^5$ due to complete exhaustion of available RAM, whereas the ML method continues to provide reliable solutions.

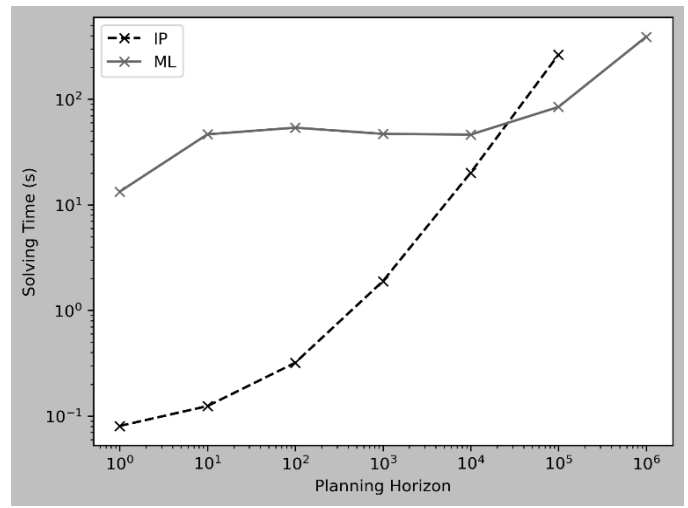

Fig. 1. Comparison of computational times of solvers

# 7. Discussion

The results demonstrate in five out of the nine scenarios examined, that the ML-based optimization method successfully reproduces the optimal solutions obtained using IP. In the remaining cases, deviations are minimal, with a maximum discrepancy of approximately 1%, while no constraint violations occur in any scenario. These findings confirm that enforcing constraints through penalty terms is sufficient to maintain feasibility, highlighting the robustness and reliability of the ML-based approach.

Table 5 exemplifies cases where the ML method does not reach exact optimality, predominantly in period 0, (except Scenario 6 where period 2 deviates). In most cases, discrepancies between IP and ML solutions are minimal, with only Scenario 2 showing a pronounced deviation.

Table 1. Vehicle Price, O&M-Cost und Utilization.

| Sce. | Optimizer | Age | | | | | | | | | | |
|---|---|---|---|---|---|---|---|---|---|---|---|---|
| | | 0 | 1 | 2 | 3 | 4 | 5 | 6 | 7 | 8 | 9 | 10 |
| 2 | IP | 0 | 0 | 0 | 11 | 0 | 0 | 0 | 0 | 0 | 1 | 0 |
| | ML | 0 | 0 | 1 | 9 | 1 | 0 | 0 | 1 | 0 | 0 | 0 |
| 4 | IP | 0 | 0 | 0 | 10 | 0 | 0 | 2 | 0 | 0 | 0 | 0 |
| | ML | 0 | 0 | 0 | 9 | 0 | 3 | 0 | 0 | 0 | 0 | 0 |
| 6 | IP | 0 | 0 | 0 | 0 | 12 | 0 | 0 | 1 | 0 | 0 | 0 |
| | ML | 0 | 0 | 0 | 0 | 13 | 0 | 0 | 0 | 0 | 0 | 0 |
| 8 | IP | 0 | 0 | 0 | 0 | 1 | 0 | 15 | 0 | 0 | 0 | 0 |
| | ML | 0 | 0 | 0 | 0 | 0 | 2 | 14 | 0 | 0 | 0 | 0 |

Although the IP solver always delivers optimal solutions, its computational performance degrades sharply as problem size grows. Additionally, the large number of branches which need to be stored leads to very high memory requirements. In contrast, the ML-based method reliably provides solutions for very large problems, as it only requires information from the current and previous optimization step, avoiding memory bottlenecks. This illustrates the clear advantage of the ML approach over IP for large-scale optimization problems.

## 8. Conclusion and Outlook

In this work, an extended IP formulation is introduced to integrate fleet renewal and upgrade options, which can be solved established solvers such as CPLEX. To overcome the extensive needs of IP w.r.t. RAM and computation time for more complex problem statements, an alternative approach based on ML is proposed. The results indicate that the IP method with a branch-and-cut solver like CPLEX is preferable only for smaller cases.

For extensive and complex product fleet management with individual upgrade strategies and additional constraints, the ML version is better scalable than the conventional solution. A machine learning approach may therefore be more suitable for vehicle fleets with numerous configurable components and variants and additional constraints.

The current contribution focuses on deterministic optimization with a predefined cost function. Future research could explore the integration of data-driven methods to better capture uncertainties and dynamic behaviors inherent in real-world systems, enabling a more nuanced and realistic modeling approach.

Furthermore, recent studies by Lee and Kim [34] and Tang et al. [35] have introduced novel ML-based techniques for solving MIP problems. A comparative analysis with these emerging methods could yield valuable insights. Moreover, the development of hybrid approaches that synergistically combine classical mathematical optimization with ML techniques holds significant potential for enhancing both scalability and solution quality.

## Declaration of generative AI and AI-assisted technologies in the writing process

During the preparation of this work the authors used DeepL and ChatGPT in order to ensure correct translations and grammar. After using these tools/services, the authors reviewed and edited the content as needed and take full responsibility for the content of the publication.

## References

[1] Khan MA, West S, Wuest T. Midlife upgrade of capital equipment: A servitization-enabled, value-adding alternative to traditional equipment replacement strategies. CIRP J Manuf Sci Technol 2020;29:232–44. https://doi.org/https://doi.org/10.1016/j.cirpj.2019.09.001.
[2] Argilovski A, Vasileska E, Jovanoski B. Enhancing manufacturing efficiency - A Lean Industry 4.0 approach to retrofitting. Mechanical Engineering-Scientific Journal 2023;41:123–9. https://doi.org/10.55302/MESJ234126672123a.
[3] Ansaripoor AH, Oliveira FS, Liret A. Recursive expected conditional value at risk in the fleet renewal problem with alternative fuel vehicles. Transp Res Part C Emerg Technol 2016;65:156–71. https://doi.org/https://doi.org/10.1016/j.trc.2015.12.010.
[4] Hopp WJ, Nair SK. Timing replacement decisions under discontinuous technological change. Naval Research Logistics (NRL) 1991;38:203–20. https://doi.org/10.1002/1520-6750(199104)38:23.0.CO.
[5] Karabakal N, Lohmann JR, Bean JC. Parallel Replacement under Capital Rationing Constraints. Manage Sci 1994;40:305–19.
[6] Hartman JC. The parallel replacement problem with demand and capital budgeting constraints. Naval Research Logistics (NRL) 2000;47:40–56. https://doi.org/10.1002/(SICI)1520-6750(200002)47:1.
[7] Chang P-T. Fuzzy strategic replacement analysis. Eur J Oper Res 2005;160:532–59. https://doi.org/https://doi.org/10.1016/j.ejor.2003.07.001.
[8] Hsu C-I, Li H-C, Liu S-M, et al. Aircraft replacement scheduling: A dynamic programming approach. Transp Res E Logist Transp Rev 2011;47:41–60. https://doi.org/https://doi.org/10.1016/j.tre.2010.07.006.
[9] Ansaripoor AH, Oliveira FS. Flexible lease contracts in the fleet replacement problem with alternative fuel vehicles: A real-options approach. Eur J Oper Res 2018;266:316–27.
[10] Foroutani A, Zarch MG, Tipaldi M, et al. A stochastic dynamic programming approach for the machine replacement problem. Eng Appl Artif Intell 2023;118:105638. https://doi.org/https://doi.org/10.1016/j.engappai.2022.105638.
[11] Parthanadee P, Buddhakulsomsiri J, Charnsethikul P. A study of replacement rules for a parallel fleet replacement problem based on user preference utilization pattern and alternative fuel considerations. Comput Ind Eng 2012;63:46–57. https://doi.org/10.1016/j.cie.2012.01.011.
[12] Büyüktahtakın İE, Hartman JC. A mixed-integer programming approach to the parallel replacement problem under technological change. Int J Prod Res 2016;54:680–95. https://doi.org/10.1080/00207543.2015.1030470.
[13] Figliozzi M, Boudart JA, Feng W. Economic and Environmental Optimization of Vehicle Fleets: Impact of Policy, Market, Utilization, and Technological Factors. Transp Res Rec 2011;2252:1–6. https://doi.org/10.3141/2252-01.
[14] Boudart J, Figliozzi M. Key Variables Affecting Decisions of Bus Replacement Age and Total Costs. Transp Res Rec 2012;2274:109–13. https://doi.org/10.3141/2274-12.
[15] Emiliano W, Alvelos F, Telhada J, et al. An optimization model for bus fleet replacement with budgetary and environmental constraints. Transportation Planning and Technology 2020;43:488–502. https://doi.org/10.1080/03081060.2020.1763656.
[16] Chung W-H, Kremer GO, Wysk RA. A dynamic programming method for product upgrade planning incorporating technology development and end-of-life decisions. Journal of Industrial and Production Engineering 2017;34:30–41. https://doi.org/10.1080/21681015.2016.1192067.
[17] Wu B, Jiang Z, Zhu S, et al. Data-Driven Decision-Making method for Functional Upgrade Remanufacturing of used products based on Multi-Life Customization Scenarios. J Clean Prod 2022;334:130238. https://doi.org/https://doi.org/10.1016/j.jclepro.2021.130238.
[18] Khan MA. Capital Equipment Upgrade Strategy in the Context of

Servitization. West Virginia University, 2020.

[19] Simons M. Comparing Industrial Cluster Cases to Define Upgrade Business Models for a Circular Economy. In: Grösser Stefan N. and Reyes-Lecuona A and GG, editor. Dynamics of Long-Life Assets: From Technology Adaptation to Upgrading the Business Model, Cham: Springer International Publishing; 2017, p. 327–56. https://doi.org/10.1007/978-3-319-45438-2_17.

[20] Albers A, Düser T, Kuebler M, et al. Upgradeable Mechatronic Systems - Definition and Model of Upgrades in the Context of the Model of SGE - System Generation Engineering. FISITA World Congress 2023 Proceedings, 2023. https://doi.org/10.46720/fwc2023-sel-009.

[21] Xing K, Abhary K. A genetic algorithm-based optimisation approach for product upgradability design. Journal of Engineering Design 2010;21:519–43.

[22] Martinez M, Xue D. A modular design approach for modeling and optimization of adaptable products considering the whole product utilization spans. Proceedings of the Institution of Mechanical Engineers, Part C 2018;232:1146–64. https://doi.org/10.1177/0954406217704007.

[23] Gadalla M, Xue D. An efficient optimisation method based on weighted AND-OR trees for concurrent reconfigurable product design and reconfiguration process planning. Int J Prod Res 2023;61:859–79. https://doi.org/10.1080/00207543.2021.2018137.

[24] Su C, Wang X. Optimizing upgrade level and preventive maintenance policy for second-hand products sold with warranty. Proceedings of the Institution of Mechanical Engineers, Part O 2014;228:518–28. https://doi.org/10.1177/1748006X14537250.

[25] Su C, Wang X. Optimal Upgrade Policy for Used Products Sold with Two-dimensional Warranty. Qual Reliab Eng Int 2016;32:2889–99. https://doi.org/https://doi.org/10.1002/qre.1973.

[26] Chung W-H, Kremer GO, Wysk RA. An Optimal Upgrade Strategy for Product Users Considering Future Uncertainty. Proceedings of the 2010 Industrial Engineering Research Conference, 2010.

[27] Medina-Oliva G, Voisin A, Monnin M, et al. Predictive diagnosis based on a fleet-wide ontology approach. Knowl Based Syst 2014;68:40–57.

[28] Gritzmann P. Grundlagen der Mathematischen Optimierung. Springer Vieweg Wiesbaden; 2013. https://doi.org/10.1007/978-3-8348-2011-2.

[29] Hartman JC. A GENERAL PROCEDURE FOR INCORPORATING ASSET UTILIZATION DECISIONS INTO REPLACEMENT ANALYSIS. Eng Econ 1999;44:217–38.

[30] Kannan Ravindran and Monma CL. On the Computational Complexity of Integer Programming Problems. In: Henn Rudolf and Korte B and OW, editor. Optimization and Operations Research, Berlin, Heidelberg: Springer Berlin Heidelberg; 1978, p. 161–72.

[31] Bengio Y, Léonard N, Courville A. Estimating or Propagating Gradients Through Stochastic Neurons for Conditional Computation 2013.

[32] Shah R, Yan M, Mozer MC, et al. Improving Discrete Optimisation Via Decoupled Straight-Through Gumbel-Softmax 2024.

[33] Bynum ML, Hackebeil GA, Hart WE, et al. Pyomo — Optimization Modeling in Python. Springer International Publishing; 2021.

[34] Kingma DP, Ba J. Adam: A Method for Stochastic Optimization 2017.

[35] Fan T-H, Chi T-C, Rudnicky AI, et al. Training Discrete Deep Generative Models via Gapped Straight-Through Estimator 2022.

[36] Rudschies W. Nio ES8: Erst mal nur in Norwegen n.d. https://www.adac.de/rund-ums-fahrzeug/autokatalog/marken-modelle/nio/nio-es8/ (accessed March 30, 2025).

[37] Lin C. The Chinese EV Company That Made Battery Swapping Work n.d. https://knowledge.insead.edu/strategy/chinese-ev-company-made-battery-swapping-work (accessed March 30, 2025).

[38] Wietschel M, Moll C, Oberle S, et al. Klimabilanz, Kosten und Potenziale verschiedener Kraftstoffarten und Antriebssysteme für Pkw und Lkw. Endbericht 2019. https://doi.org/10.24406/publica-fhg-299856.

[39] Mendelevitch R, Repenning J, Matthes FChr, et al. Treibhausgas-Projektionen 2024 für Deutschland - Rahmendaten. Dessau-Roßlau: Umweltbundesamt; 2024.

[40] Wang B, Fu J, Zhang H, et al. Closing the Gap Between the Upper Bound and the Lower Bound of Adam's Iteration Complexity 2023.

[41] Lee T-H, Kim M-S. RL-SPH: Learning to Achieve Feasible Solutions for Integer Linear Programs 2025.

[42] Tang B, Khalil EB, Drgoňa J. Learning to Optimize for Mixed-Integer Non-linear Programming with Feasibility Guarantees 2025.